\documentstyle{amsppt}
\hoffset =.1in
\NoBlackBoxes
\magnification=\magstep 1
\nologo

\baselineskip18pt
\pageheight{25truecm}
\pagewidth{15truecm}
\vcorrection{-1truecm}

\nopagenumbers

\define\Proj{\operatorname{Proj}}
\define\lin{\operatorname{lin}}



\topmatter
\title  In the Amemiya-Ando problem 3 is enough
\endtitle
\author Adam Paszkiewicz
\endauthor

\endtopmatter

\footnote""{{\it 2010} AMS subject classification: 47H09.}
\footnote""{Research supported by  grant N N201 605840.}

\flushpar{\bf 1. The main result}.
We shall prove the following 

\proclaim{1.1. Theorem} In any infinite dimensional Hilbert space $\Cal H$ there exist orthogonal projections $Q_1, Q_2, Q_3$ such that a sequence $P_n\dots P_1x$ diverges in norm for some $P_n\in\{Q_1,Q_2, Q_3\}$, $n\geq 1$,  $x\in \Cal H.$  

\endproclaim

Moreover, quait simply sequence $(P_n)$ can be taken. Namely there exist $Q_1, Q_2, Q_3\in \Proj \Cal H$, $x\in \Cal H$ and numbers $m(i)$, $p(i) \geq 1$ such that
$$\prod_{1\leq i\leq j} (Q_3(Q_2Q_1Q_2)^{p_i})Q_3)^{m_i} x$$
diverges in norm for $j\to \infty$. In the whole paper $\prod_{1\leq i\leq j} A_i$ will be a short notation for $A_j\dots A_1$, for any $A_1,\dots, A_j\in B(\Cal H).$

By capital letters $D,  E, F, G, P, Q, R$ we shall always denote orthogonal projections in $\Cal H$. By small letters $e, f, p, q$ we shall always denote some vectors in $\Cal H$, of norm $1$; $e'$ is always orthogonal to $e$ ($f'\bot f$ and so on). Then $\hat e\cdot= \langle \cdot, e\rangle e$ is a one-dimensional projection and always $\hat e'\bot \hat e$, then $\hat f\leq \hat e +\hat e'$ means that $f\in \lin (e, e').$

Theorem 1.1 it is a consequence of the following two lemmas. 

\proclaim{1.2. Lemma (The main difficulty)} Let $E =\sum\hat e_i$, $\dim E^\bot =\infty$ and let $\epsilon_i >0$, $i\geq 1$. There exist numbers $m_{i s}\geq 1$ and projections $P_{is}$, $1\leq s\leq \sigma_i$, satisfying
$$||\prod_{1\leq s\leq \sigma_i} (E P_{i s}E)^{m_{is}} e_i -e_{i+1}|| <\epsilon_i, \tag1$$
and all projections $\{P_{is}; 1\leq s\leq \sigma_i$, $i\geq 1\}$ can be ordered into a decreasing sequence $P_1\geq P_2 \geq\dots$, and $\dim P_1^\bot =\infty$, $\dim(P_m -P_{m+1}) <\infty$, $m\geq 1.$

\endproclaim

\proclaim{1.3. Lemma (cf. Lemma 2 in [1])} For any projections $P_1\geq P_2\geq \dots$, $\dim P_1^\bot =\infty$, $\dim(P_m -P_{m+1}) <\infty$,  and numbers $\delta_1, \delta_2,\dots >0$, there exist projection $Q$ and numbers $p_1, p_2,\dots \geq 1$ satisfying 
$$||(P_1 QP_1)^{p_k} - P_k|| < \delta_k, \quad k\geq 1.$$
\endproclaim

\flushpar{\bf 1.4. Proof of Theorem 1.1}

It is rather immediate consequence of Lemmas 1.2, 1.3. Projections $P_{is}$ satisfying (1) can be found for $\epsilon_i =1/4\cdot 2^i$ and one can assume that $\dim P_1^\bot =\infty$ for the maximal element $P_1 = \bigvee P_{is}$. By Lemma 1.3 there exists projection $Q$ satisfying $||P_{is} - (P_1 Q P_1)^{p_{is}}|| <\delta_{is}$ for arbitrary small $\delta_{is} >0$, $1\leq s \leq s_i$. Thus the estimation
$$||\prod_{1\leq s\leq \delta_i} (E P_{is} E)^{m_{is}} - A_i|| < 1/4 \cdot 2^i$$
can be obtained for
$$A_i = \prod_{1\leq s\leq\sigma_i} (E(P_1 Q P_1)^{p_{is}} E)^{m_{is}}$$
if only we take $\delta_{is}$ small enough. In consequence, a relation
$$||A_{i-1}\dots A_1 e_1 - e_i|| \leq \frac12 - \frac{1}{2^i}\tag2$$
implies $||A_i\dots A_1 e_1 - e_{i+1}|| =||A_i(A_{i-1}\dots A_1 e_1 - e_i) + (A_i -\prod_{1\leq s \leq\sigma_1} (E P_{is} E)^{m_{is}})e_i +
\prod_{1\leq s \leq \sigma_i} (E P_{is} E)^{m_{is}} e_i - e_{i+1}|| < (\frac12 -\frac{1}{2^i}) +\frac{1}{4\cdot 2^i} +\frac{1}{4\cdot 2^i} =\frac12 -\frac{1}{2^{i+1}}$. The estimation (2) is also valid for $i =1$ and thus for any $i\geq 1$. It means that $x= e_1$ and $\{Q_1, Q_2, Q_3\} =\{E, P_1, Q\}$ an be taken.

\bigskip

\flushpar{\bf 2. The overcome of the main difficulty}
\medskip

The proof of Lemma 1.2 is more delicate. At first we show that it can be obtained from the following lemma (being a stronger version of Lemma 1 in [1]):

\proclaim{2.1. Lemma} For any $\epsilon >0$ there exists $t(\epsilon) >0$ such that for any projections
$$E + G = Q_0\geq \dots \geq Q_{t(\epsilon)}= E$$
satisfying $\dim(Q_{t-1} - Q_t) =1$, $1\leq t\leq t(\epsilon))$, $E =\hat e+\hat e'$, and for $\eta >0$ there exist a unitary operator $V$ and numbers $n(1),\dots, n(t(\epsilon))$ satisfying
$$V Q_0^\bot = Q_0^\bot,$$
$$||V - {\bold 1}|| < \eta,$$
$$||\prod_{1\leq t\leq t(\epsilon)} (E V Q_t V^* E)^{m(t)}e - e'|| <\epsilon.$$

\endproclaim

Passing from Lemma 2.1 to Lemma 1.2 can be simplified by a formal use of the following two lemmas. We shall need, in particular some (rather special) assumptions which quaranties the estimations
$$||(\tilde V\tilde U\tilde P\tilde U^* \tilde V^* - U P U^*) P_0|| < \eta,\tag3$$
$$(\tilde V\tilde U \tilde Q\tilde U^* \tilde V^* - U Q V)Q_0 =0\tag4$$
for projections $P, P_0, \tilde P, Q, Q_0, \tilde Q$ and unitary $U, \tilde U, V, \tilde V.$

\proclaim{2.2. Lemma} A. If $\tilde U\tilde Q\tilde U^* =\tilde Q$, $\tilde V Q_0 = V Q_0 = Q_0 \tilde V$, and $Q_0 \tilde Q = Q$, then (4).

B. If $\tilde U P_0 = U P_0 = P_0\tilde U$, $P_0\tilde P = P$, and $||\tilde V P_0 - P_0||\leq \eta/2$ for $\eta >0$, then (3).

\endproclaim

\demo{Proof A} As $V, V_1$ coincide on the space $Q_0 \Cal H$ and commute with $Q_0$, we have
$$\tilde V \tilde U \tilde Q \tilde U^* \tilde V^* Q_0 = \tilde V \tilde Q \tilde V^* Q_0 = \tilde V \tilde Q Q_0 \tilde V^* Q_0 = V_1 Q V_1^* Q_0.$$

B. By properties of $U, U_1, P_0$ we have (analogically) $\tilde U \tilde P \tilde U^* P_0 = U P U^* P_0$. As $||V P_0 - P_0|| < \eta/2$, we have also $||V^* P_0 - P_0|| = ||V^*(P_0 - V P_0)|| <\eta/2$ and, for $\tilde p = \tilde U \tilde P \tilde U^*$, $p = U P U^*$,
$$\tilde p P_0 = P_0 \tilde p P_0 = p,$$
$$||\tilde V \tilde p \tilde V^* P_0 - p|| = ||\tilde V \tilde p(\tilde V^* P_0 - P_0) + (\tilde V P_0 - P_0)\tilde p P_0 +(P_0 \tilde p P_0 - p)|| <\eta/2 +\eta/2 +0.$$

\enddemo

\proclaim{2.3. Lemma} Let $A_1,\dots, A_M$ be a system of operators in $\Cal H$, $||A_m||\leq 1$.

A) If $(A_m - B_m)Q_0 =0$, $A_mQ_0 =Q_0 A_m Q_0$, $B_m\in B(\Cal H)$ for $1\leq m\leq M$, $Q_0\in\Proj \Cal H$ then
$$(\prod_{1\leq m\leq M} B_m)Q_0 = (\prod_{1\leq m\leq M} A_m)Q_0.\tag5$$

B) For any $\epsilon >0$ there exists $\gamma >0$ such that $||(A_m - B_m) P_0|| <\gamma$, $A_m P_0 = P_0 A_m P_0$, $||B_m|| \leq 1$, $B_m\in B(\Cal H)$ for $1\leq m\leq M$, $P_0\in\Proj \Cal H$ imply
$$||(\prod_{1\leq m\leq M} A_m) P_0 - (\prod_{1\leq m\leq M} B_m) P_0|| <\epsilon.\tag6$$

\endproclaim

\demo{Proof} A) We have $A_{M-1}\dots A_1 Q_0 = Q_0 A_{M-1}\dots A_1 Q_1$, $B_{M-1}\dots B_1 Q_0 =$ $Q_0 B_{M-1}\dots$ $B_1 Q_0$. Thus $A_{M-1}\dots A_1 Q_0 =$ $B_{M-1}$ $\dots B_1 Q_0$ implies (5), and thus (5) is proved by induction.

B) If the condition B) is valid for $M-1$ instead of $M$ then for $\epsilon >0$ there exists $\gamma >0$ such that $\gamma <\epsilon /2$ and $||A - B|| <\frac{\epsilon}{2}$ for
$$A = (\prod_{1\leq m\leq M-1} A_m) P_0,$$
$$B = (\prod_{1\leq m\leq M-1} B_m) P_0$$
(for $A_m$, $B_m$ satisfying suitable assumptions for $M-1$ instead of $M$).

Then $||(A_M - B_M)P_0|| <\gamma$, $||B_M|| \leq 1$ imply $||A_M A-B_M B||\leq ||A_M A - A_M P_0 A|| + ||A_M P_0A - B_M P_0 A|| + ||B_MP_0 A-B_MA|| + ||B_M A-B_M B|| < 0 +\gamma +0 +\epsilon/2 <\epsilon$. The condition B) is proved by induction.

\enddemo

\flushpar{\bf 2.4. Proof of Lemma 1.2.}

Let us take numbers $s_k = t(\frac12 \epsilon_{2k-1})$, $t_k = t(\epsilon_{2k})$, $k\geq 1$, according to Lemma 2.1, and mutually orthogonal projections $E; F_1, F_2,\dots$, $G_1, G_2,\dots$ satisfying
$$\dim(E +\sum(F_k +G_k))^\bot =\infty,$$
$$E = \sum_{i\geq 1} \hat e_i,\qquad (e_i) \text{-- O.N. system},$$
and such that
$$\dim(P_{s-1}^k - P_s^k) =\dim(Q_{t-1}^k - Q_t^k) =1\tag13$$
for some projections
$$\aligned D_k +F_k = P_0^k \geq \dots \geq P_{s_k}^k = D_k, &\quad D_k = \hat e_{2k-1} +\hat e_{2k},\\
E_k + G_k = Q_0^k \geq \dots \geq Q_{t_k}^k = E_k, &\quad E_k = \hat e_{2k} +\hat e_{2k+1}.\endaligned\tag14$$
Then for fixed $k\geq 1$ and $\eta_k > 0$ (defined later), one can find a unitary operator $V_k$ and numbers $n(k, t) \geq 1$ satisfying
$$||\prod_{1\leq t\leq t_k} (E_k(V_k Q_t^k V_k^*) E_k)^{n(k, t)}e_{2k} - e_{2k+1}|| <\epsilon_{2k},\tag7$$
$$\aligned & V_k Q_0^{k\bot} = Q_0^{k\bot},\\
& ||V_k -{\bold 1}|| <\eta_k.\endaligned\tag12$$
Analogically, one can find a unitary $U_k$ and numbers $m(k, s)\geq 1$ satisfying
$$||\prod_{1\leq s\leq s_k} (D_k(U_k P_s^k U_k^*) D_k)^{m(k, s)}e_{2k-1} - e_{2k}|| < \frac12 \epsilon_{2k-1},\tag8$$
$$U_k P_0^{k\bot} = P_0^{k \bot}.$$
In particular, we have well defined unitary operators
$$\tilde U = \prod_{k\geq 1} U_k,\qquad \tilde V =\prod_{k\geq 1} V_k,$$
and
$$\tilde U P_0^k = U_k P_0^k = P_0^k \tilde U, \quad \tilde V Q_0 = V_k Q_0^k = Q_0^k\tilde V,$$
$$||\tilde V P_0^k - P_0^k|| \leq \eta_{k-1} \vee \eta_k, \ k\geq 1  \ \text{(assuming that $\eta_0 =0$)}.$$
The last estimation is a consequence of $||\tilde V(Q_0^k \vee Q_0^{k-1}) - Q_0^k \vee Q_0^{k-1}||\leq \eta_{k-1} \vee \eta_k$, $P_0^k\leq Q_0^{k} \vee Q_0^{k-1} +F_k$, $k\geq 1$, where we put $Q_0^0 =0$. Moreover,
$\tilde V F_k = F_k$ (as $F_k\leq Q_0^{l \bot}$ and $V_l F_k = F_k$ for $l, k\geq 1$).

Now we pass to a system of projections, which can be ordered into a decreasing sequence:
$$\tilde P_s^k = P_s^k +\sum_{l\geq k+1} P_0^l,$$
$$\tilde Q_t^k = P_0^k + Q_t^k - E_k +\sum_{l\geq k+1} P_0^l.$$
Indeed we have
$$\multline \tilde Q_0^1 \geq \dots \geq Q^1_{t_1}\geq \tilde P_0^1 \geq\dots \geq \tilde P^1_{s_1}\geq\dots\\
\dots \geq \tilde Q_0^k \geq \tilde Q^k_{t_k} \geq \tilde P_0^k \geq \dots \geq \tilde P^k_{s_k}\geq \dots\endmultline\tag11$$
Then
$$\tilde U \tilde Q_t^k \tilde U^*= \tilde Q_t^k,$$
$$Q_0^k \tilde Q_t^k = Q_t^k,\qquad P_0^k\tilde P_s^k = P_s^k.$$

Lemma 2.2 A. gives
$$(\hat Q_t^k - V_k Q_t^k V^*_k) Q_0^k =0,$$
for $\hat Q_t^k = \tilde V \tilde U \tilde Q_t^k \tilde U^*\tilde V^*$, and Lemma 2.2 B. gives
$$||(\hat P_s^k - U_k P_s^k U_k^*) P_0^k|| < 2(\eta_{k-1} \vee \eta_k)$$
for $\hat P_s^k = \tilde V \tilde U \tilde P_s^k \tilde U^* \tilde V^*.$

Now we use Lemma 2.3 A) with $A_m$ being $V_k Q_t^k V_k^*$ or $E_k$, $B_m$ being $\hat Q_t^k$ or $E$ (respectively), and this gives
$$(\prod_{1\leq t\leq t_k} (E \hat Q_t^k E)^{n(k, t)}) Q_0^k = (\prod_{1\leq t\leq t_k} (E_k V_k Q_t^k V_k^* E_k)^{n(k, t)})Q_0^k.\tag9$$
Analogically, Lemma 2.3 B), with $A_m$  being $U_k P_s^k U_k^*$ or $D_k$, $B_m$ being $\hat P_s^k$ or $E$ gives
$$\aligned &||(\prod_{1\leq s \leq s_k} (E \hat P_s^k E)^{m(k, s)})P_0^k \\
- &(\prod_{1\leq s\leq s_k} (D_k U_k P_s^k U_k^* D_k)^{m(k, s)}) P_0^k|| <\epsilon_{2k-1}/2,\endaligned\tag10$$
if only $\eta_{k-1}\vee \eta_k \leq \gamma_k$ for suitable small $\gamma_k >0$.

Now it is a good moment to define numbers $\eta_k$, such that all relations (7), (8), (9), (10) can be satisfied. We recall that, for a sequence $\epsilon_i$, $i\geq 1$, the numbers $s_k$, $t_k$ are given by Lemma 2.1. Then operators $U_k$ and numbers $m(k, s)$, $1\leq s\leq s_k$, satisfying (8) can be immediately found (for any system $P_s^k$).

Then we take $\gamma_k=\gamma$, given by Lemma 2.3 B) for $\epsilon =\epsilon_{2k-1}/2$, and for $M$ being a number of terms in the product $\prod_{1\leq s\leq s_k} (D_k R_s^k D_k)^{m(k, s)}$ (with $R_s^k = U_k P_s^k U_k^*$). We put $\eta_k =\gamma_k \wedge \gamma_{k+1}$, $k\geq 1$, and take operators $V_k$ satisfying (7), (12), given by Lemma 2.1.

The relation (9) can be obtained immediately and (10) can be also obtain as a consequence of $\gamma_k =\eta_k\vee \eta_{k-1}$, $k\geq 1$ ($\eta_0 =0$).

The formula (1) is a consequence of (7), (8) and (9), (10), after the following change of notations
$$P_{2k-1, s} = \hat P_s^k,\qquad 1\leq s \leq \sigma_{2k-1} := s_k,$$
$$P_{2k, t} = \hat Q_t^k,\qquad 1\leq t \leq \sigma_{2k} := t_k.$$
Moreover by (11), we have a decreasing ordering
$$\multline(P_m) = (P_{2, 0}, \dots, P_{2, \sigma_2}; P_{1, 0}, \dots,  P_{1, \sigma_1}; \dots\\
\dots; P_{2k, 0}, \dots, P_{2k, \sigma_{2k}}; P_{2k-1, 0}, \dots, P_{2k-1, \sigma_{2k-1}}; \dots)\endmultline$$
and $\dim(P_m - P_{m+1}) =1$ or $2$, by (13), (14). Using once more a definition of $\tilde U, \tilde V$ we have
$$\hat Q_0^2 =\tilde V \tilde U Q_0^2 \tilde U^* \tilde V^* =Q_0^2,$$
because $Q_0^2 = E +\sum(F_k +G_k),$ and
$$\dim P_{2, 0}^\bot =\dim \hat Q_0^{2 \bot} =\infty.$$
The proof is finished.

\bigskip

\flushpar{\bf Referencers}

\ref\key 1
\by A. Paszkiewicz
\paper The Amemiya-Ando conjecture falls
\jour preprint
\endref

\bye